\documentclass[a4paper,12pt,oneside]{amsart}

\usepackage{amsmath,amsfonts,amssymb,amsthm,amsopn,amsxtra}

\newtheorem*{Thm}{Theorem}    
\newtheorem{Lem}{Lemma}
\newtheorem{Pro}{Proposition}
\newtheorem*{Cor}{Corollary}

\theoremstyle{remark}

\DeclareMathOperator*{\esup}{ess\,sup}
\DeclareMathOperator{\ch}{ch}

\DeclareMathOperator{\tah}{th}
\DeclareMathOperator{\tr}{tr}
\DeclareMathOperator{\im}{im}
\DeclareMathOperator{\re}{Re}
\DeclareMathOperator{\sgn}{sgn}

\newcommand{\Cal}{\mathcal}
\newcommand{\fr}{\mathfrak}
\newcommand{\R}{\mathbb{R}}        
\newcommand{\C}{\mathbb{C}}        
\newcommand{\Z}{\mathbb{Z}}        
\newcommand{\N}{\mathbb{N}}
\newcommand{\VN}{\mathit{VN}}
\newcommand{\SL}{\mathit{SL}}
\newcommand{\PSL}{\mathit{PSL}}

\begin{document}
\title{On multipliers and completely bounded multipliers -- the case
$\SL(2,\R)$}
\author{Viktor  Losert}
\address{Institut f\"ur Mathematik, Universit\"at Wien, Strudlhofg.\ 4,
  A 1090 Wien, Austria}
\email{Viktor.Losert@UNIVIE.AC.AT}
\thanks{Summary of talks given at the School of Mathematics, University of
Leeds, 28 May--\;2 June 2010; revised version.}
\date{ December 2014}


\maketitle

\baselineskip=1.3\normalbaselineskip
\begin{list}{}{
\itemindent-5mm \labelwidth0pt \setlength{\leftmargin}{5mm} }
\item $A(G)$ Fourier algebra of a locally compact group $G$\,. $B(G)$
Fourier-Stieltjes algebra.
\item $A(G)''$ bidual of $A(G)$ with (first) Arens product $\odot$\;.
\item $M(A(G))$ multipliers of $A(G)$ with norm $\lVert\;\rVert_M$\,.
Every $f\in M(A(G))$ is given by (and identified with) a bounded continuous
function on $G$. It extends to $A(G)''$ and this is again denoted by
$f\odot\xi$ for $\xi\in A(G)''$ (bidual mapping).
\item $M_0(A(G))$ completely bounded multipliers of $A(G)$ with norm
$\lVert\;\rVert_{M_0}$ \,(see \cite{CH} for basic properties).
\item $\VN(G)$ group von Neumann algebra (generated by the left regular
representation on $L^2(G)$\,), we use the standard identification with
the dual space $A(G)'$.
\item $C_0(G)$ continuous functions on $G$ vanishing at infinity.
\item $\Cal B(\Cal H)$ bounded linear operators on a Hilbert space $\Cal H$\,.
\item $\Cal N(\Cal H)$ nuclear operators (trace class), identified with the
predual $\Cal B(\Cal H)_*$\,, using $(t,s)=\tr(t\,s)$.
\end{list}

\vskip 1mm\noindent
For $G=\SL(2,\R)$ (real 2x2-matrices of determinant one), let $K$ be
the subgroup of rotations
$k_{\varphi}=
\begin{pmatrix}\cos\varphi & -\sin\varphi\\\sin\varphi& \cos\varphi
\end{pmatrix}$\vspace{.5mm} and $H$ the subgroup of matrices
$\begin{pmatrix}a & 0\\b& \frac1a
\end{pmatrix}$
with $a>0,\ b\in\R$\,. Recall (part of the Iwasawa decomposition) that
$G=KH$\,, the decomposition of the elements $x=kh$ being unique. We
parametrize the dual group $\widehat K$ of the compact abelian group $K$
by $\chi_j(k_{\varphi})=e^{ij\varphi}$  ($j\in\Z,\ \varphi\in\R$).
For a bounded continuous function $f$ on $G$, $m,n\in\Z$ put
$f_{mn}=(\chi_n*f*\chi_m)\,\vert H$\,.
\\
To simplify, we describe the main result in the case of
$\PSL(2,\R)=\SL(2,\R)/\{\pm I\}$
(projective special linear group; $\{\pm I\}$ being the centre of
$\SL(2,\R)$\,).  Then $K$ is replaced by  $K/\{\pm I\}$ and $f_{mn}$ is
defined only for even $m,n$\,.
\begin{Thm}
For $G=\PSL(2,\R)$ we have $M(A(G))=M_0(A(G))$.
\[
\lVert f\rVert_M=\lVert f\rVert_{M_0}\qquad
\text{holds for all \ }f\in M(A(G)).
\]
Then, putting $\Cal A=\VN(H)\bar\otimes \Cal B(l^2(2\Z))$,
the following statements are equivalent for $f\in C_0(G)$:\vspace{-1mm}
\begin{enumerate}
\item  $f\in M(A(G))$
\item
$(f_{mn})_{m,n\in2\Z}$ defines an element of the predual $\Cal A_*$ of
$\Cal A$ and \\ $\theta_f(e_{kl})=(f_{m-k\,n-l})_{m,n\in2\Z}$ extends to a
continuous linear mapping \\ $\Cal N(l^2(2\Z))\to\Cal A_*$\,.
\end{enumerate}
Furthermore\quad $\lVert f\rVert_M=\lVert f\rVert_{M_0}=
\lVert \theta_f\rVert$ \ holds.
\\
For general $f\in M(A(G))$, we have that $\lambda=\lim_{x\to\infty}f(x)$
exists. Then\linebreak $f-\lambda\in M(A(G))\cap C_0(G)$ and
$\lVert f\rVert_M=\lVert f-\lambda\rVert_M+\lvert\lambda\rvert$\,.
\end{Thm}
\noindent
As in \cite{Ta}\,p.184, elements $t\in \VN(H)\bar\otimes \Cal B(l^2(\Z))$
are described by matrices\linebreak
$(t_{mn})_{m,n\in\Z}$\,, where $t_{mn}\in \VN(H)$
and similarly $s\in\Cal A_*$ is given by  $(s_{mn})_{m,n\in2\Z}$\,, where
$s_{mn}\in A(H)$. This amounts to
$ (t,s)=\sum_{m,n}(t_{mn}\,,s_{nm})$ \ (compare \cite{Ta}\linebreak p.65(18)).
Analogously
for elements of $\Cal B(l^2(\Z))$ and $\Cal N(l^2(\Z))$. $e_{kl}$ denotes
the element of $\Cal N(l^2(\Z))$ given by the matrix with $1$ at $(k,l)$ and
$0$ elsewhere.\vspace{1mm}
 
For $G=\SL(2,\R)$ one has to observe that $f_{mn}=0$ whenever $m-n$ is odd.
With $\theta_f(e_{kl})=(f_{m-k\,n-l})_{m,n\in\Z}$ for $k,l\in2\Z$ and
$\Cal A=\VN(H)\bar\otimes \Cal B(l^2(\Z))$ (or the subalgebra of operators
commuting with the projection $l^2(\Z)\to l^2(2\Z)$),
the Theorem holds in the same way.
Similarly, the Theorem holds for all connected groups $G$ that are locally
isomorphic to $\SL(2,\R)$ and have finite centre. With some modifications,
one can find presumably also a version for the universal covering group of
$\SL(2,\R)$.

For general $G$\,, we have
$A(G)\subseteq B(G)\subseteq M_0(A(G))\subseteq M(A(G))$. When $G$ is
amenable (e.g. abelian or compact), $M(A(G))=B(G)$ holds.
When $G$ is non-amenable (e.g., $\SL(2,\R)$ or the discrete free group $F_2$),
it is known that $B(G)$ is a proper subspace of $M_0(A(G))$. For a general
discrete group $G$, containing $F_2$ as a subgroup, Bozejko (1981) has
shown that $M_0(A(G))$ is a proper subspace of $M(A(G))$.

If $K$ is a compact subgroup of some locally compact group $G$, a function $f$
on $G$ is called {\it radial} (with respect to $K$) or $K$--bi-invariant, if
$f(k_1xk_2)=f(x)$ holds for all $x\in G,\ k_1,k_2\in K$\,. If there
exists a closed amenable subgroup $H$ of $G$ such that $G=KH$ holds
set-theoretically, then for a radial
function $f$\,, Cowling and Haagerup \cite{CH} have shown that the
following conditions are equivalent:
\[
\text{{\rm (i)} \ $f\in M(A(G))$\qquad
{\rm (ii)} \ $f\in M_0(A(G))$\qquad
{\rm (iii)} \ $f\vert H\in B(H)$}
\]
(with equality of norms).
This applies, in particular, for a semisimple Lie group $G$ with finite centre,
$K$ a maximal compact subgroup.

For $G=\SL(2,\R)$ (or $\PSL(2,\R)$) and $m,n\in\Z$, using our notation above,
we call $f\ \ (m,n)$-radial,
if $f(k_1xk_2)=\chi_m(k_1)f(x)\chi_n(k_2)$ holds for all
$x\in G,\ k_1,k_2\in K$\,. Then the same equivalence as above holds for
$(m,n)$-radial functions $f$ and for $(m,n)\ne(0,0)$ one even gets
(by our Theorem) \,$f\vert H\in A(H)$. Furthermore, one can show that
the closure of $A(G)$ in $M(A(G))\cap C_0(G)$ contains all $K$-finite
functions (i.e., all $f$ for which $f_{mn}=0$ apart of finitely many
$(m,n)$\,). In some cases (e.g., when $f$ is of diagonal type, i.e.,
$f_{mn}=0$ for $m\neq n$), one can show that the condition
$(f_{mn})_{m,n\in2\Z}\in\Cal A_*$ of the Theorem is already sufficient to
conclude that $f\in M(A(G))$. But one can show that there are $f\in C_0(G)$
satisfying $(f_{mn})_{m,n\in2\Z}\in\Cal A_*$ but $f\notin M(A(G))$, i.e.,
the first condition of  \thetag{2} above is not sufficient in general
(contrary to the assertion in the first version of this draft). Observe that
our definition of $f_{mn}$ involves transposition, i.e., $f_{mn}$ is
(the restriction of) an $(n,m)$-radial function.

\medskip\setlength{\parindent}{0pt}
On the following pages, we indicate the {\sc proof} of the Theorem:
\\[1mm]
In one direction, we use a slight extension of the results of
\cite{S}\,Th.\,2.1 and 3.1 (the proofs there work quite similarly).
\begin{Pro}\label{M0}
Let $\Cal A_0, \Cal B_0$ be unital C*-algebras, $\Cal H_0, \Cal H$ be Hilbert
spaces and assume that given are *-representations of $\Cal A_0$ and $\Cal B_0$
on both spaces $\Cal H_0$ and $\Cal H$ (the operators defined by elements
of $\Cal A_0$ and $\Cal B_0$ will be denoted by the same letters). Let
$\Cal E$ be a linear subspace of $\Cal B(\Cal H_0)$ such that $aeb\in\Cal E$
for $a\in\Cal A_0,\,e\in\Cal E,\,b\in\Cal B_0$ (i.e., $\Cal E$ is an
$\Cal A_0$\,-\,$\Cal B_0$ submodule of $\Cal B(\Cal H_0)$). Let
$\phi\!:\Cal E\to\Cal B(\Cal H)$ be a linear map satisfying
$\phi(aeb)=a\phi(e)b$ for $a\in\Cal A_0,\,e\in\Cal E,\,b\in\Cal B_0$\,.
\item[(i)]\;If $\phi$ is bounded and there exist $\xi,\eta\in\Cal H$ such
that $\Cal A_0\xi$ and $\Cal B_0\eta$ are dense in $\Cal H$\,, then
$\phi$ is completely bounded and $\lVert\phi\rVert_{cb}=\lVert\phi\rVert$.
\item[(ii)]\;If $\Cal E=\Cal K(\Cal H_0)$, $\phi$ is completely bounded,
then (for some index set $I$) there exist families
$(s_i),(t_i)\subseteq\Cal B(\Cal H,\Cal H_0)$ such that
$bs_i=s_ib,\ at_i=t_ia$ for $a\in\Cal A_0,\,b\in\Cal B_0,\,i\in I$ (i.e.,
$s_i\,,t_i$ are intertwining operators for the actions of $\Cal A_0$ and
$\Cal B_0$),
$\sum_i s_i^*s_i\,,\sum_it_i^*t_i\in\Cal B(\Cal H)$,
$\lVert\sum_i s_i^*s_i\rVert\;\lVert\sum_it_i^*t_i\rVert=
\lVert\phi\rVert_{cb}^2$ and $\phi(k)=\sum_i t_i^*ks_i$ holds for all
$k\in\Cal K(\Cal H_0)$.
\end{Pro}\noindent

In fact, we will use statement (ii) slightly more generally for 
$\Cal E=\Cal K(\Cal H_1)\oplus\Cal K(\Cal H_2)$ where $\Cal H_1,\Cal H_2$
are $\Cal A_0,\Cal B_0$-invariant subspaces of
$\Cal H_0=\Cal H_1\oplus\Cal H_2$\,. As a further
extension (for the case of the universal covering group of $\SL(2,\R)$),
this holds when $\Cal E$ is a von Neumann subalgebra of $\Cal B(\Cal H_0)$,
$\phi$ is w*-continuous and the operators on $\Cal H_0$ defined by
$\Cal A_0,\Cal B_0$ belong to $\Cal E$\,.
\\
To show that $(2)\Rightarrow(1)$ assume that $f\!:G\to\C$ is continuous,
$(f_{mn}\vert H)_{m,n\in2\Z}$ defines an element of the predual of
\,$\VN(H)\,\bar\otimes\,\Cal B(l^2(2\Z))$ and $\theta_f$ is defined as in
$(2)$, with $\lVert\theta_f\rVert=c$\,. As 
explained later (after Lemma 6), $\VN(H)$ is isomorphic (as a W*-algebra) to
$\Cal B\bigl(L^2(]-\infty,0])\bigr)\oplus\Cal B\bigl(L^2([0,\infty[)\bigr)$
using a certain representation $\pi_0$ of $H$ on $L^2(\R)$.
We apply now Proposition 1, taking
$\Cal H=l^2(2\Z)$, $\Cal H_0=L^2(\R)\otimes l^2(2\Z)$,
$\Cal A_0=\Cal B_0=C^*(2\Z)$ (operating by translations)
and for $\phi$ (the restriction of) the dual mapping of $\theta_f$\,.
Obtaining $(s_i),(t_i)$ as above, we put $p_i=s_i(e_0), q_i=t_i(e_0)$
($e_n$ denoting the standard basis of $l^2(2\Z)$). Then
$p_i=(p_{ik})_{k\in2\Z}\,,\,q_i=(q_{ik})_{k\in2\Z}\in
L^2(\R)\otimes l^2(2\Z)$.
This gives $f_{mn}(h)=\sum_i(\pi_0(h)\,p_{im}\!\mid q_{in})$ for $h\in H$\,.
For $x=hz\in G$, where $h\in H,\,z\in K$ put
$P(x)=(\sum_k\chi_k(z)\pi_0(h)p_{ik})_{i\in I}\,,\,
Q(x)=(\sum_k\chi_k(z)\pi_0(h)q_{ik})_{i\in I}$.
Then (observe that $\sum_i s_i^*s_i\,,\sum_it_i^*t_i\in \VN(2\Z)$) it
follows that $P,Q$ define (a.e.) bounded measurable functions
$G\to L^2(\R)\otimes l^2(I)$,
$\esup\limits_{x,y\in G}\,\lVert P(x)\rVert\ \lVert Q(y)\rVert=c$ and
$f(y^{-1}x)=(\,P(x)\mid Q(y)\,)$ \;holds a.e.\;on $G\times G$. By
\cite{CH}\,p.\,508, we get $f\in M_0(A(G))$ and $\lVert f\rVert_{M_0}\le c$
 \;(to avoid problems of convergence and sets of measure zero, one can use
Fejer sums and first consider the $K$-finite case [i.e. where only
finitely many $f_{mn}$ are non-zero]~). One can also show that $f\in C_0(G)$.

\medskip
For the other direction, we start by recalling the description of the
irreducible unitary {\bf representations} (going back to Bargmann).
We use (essentially) the notations (and
parametrization) of Vilenkin \cite{V}.

Put $\Cal H=L^2(\R)$ (for ordinary Lebesgue measure),
$g=\begin{pmatrix}\alpha&\beta \\ \gamma& \delta
\end{pmatrix}$,
\[
\bigl(T_l(g)f\bigr)\,(x)=
f\Bigl(\frac{\alpha x+\gamma}{\beta x+\delta}\Bigr)\,
\lvert\beta x+\delta\rvert^{2l} \qquad\text{for}\quad f\in\Cal H\,.
\]
For $l=-\frac12+i\lambda$ with $\lambda\in\R$ this gives unitary (strongly
continuous, irreducible) representations of $\SL(2,\R)$
(first {\it principal series}). $-\frac12\pm i\lambda$ gives equivalent
representations, hence it will be enough to consider $\lambda\ge0$.
\\
For $l\in\Z$ one gets the {\it discrete series} (but here the inner
product has to be changed to make $T_l$ unitary, also restricting
to subspaces of $\Cal H$ for irreducibility; see below).
\\
Further cases for unitary representations are $l\in]-1,0[$, which gives the
complementary series (again with a different inner product) and, finally,
there is also the trivial (one-dimensional) representation. These are all
the irreducible unitary representations defined on $\PSL(2,\R)$.
\\
$T_l$ arises from the right action of $\SL(2,\R)$ on $\R^2$ (and the
corresponding action on the projective line). In the notation of \cite{V}
this is $T_{\chi}$ with $\chi=(l,0)$ (the second parameter can be used
to describe further representations of $\SL(2,\R)$ and other covering groups).
Integer case: \ for $l\ge0$\,, we take $T_l$  to be only the
part $T_{\chi}^-$ (notation of \cite{V}) and for $l<0$ the part $T_{\chi}^+$.
Thus $T_{-l-1}$ is (equivalent to) the conjugate representation of $T_l$\,.
\vspace{1mm}

Multiplication in $A(G)$ and $B(G)$ corresponds to {\bf tensor products} of
representations. For $\SL(2,\R)$ the decompositions have been determined by
Pukanszky (1961). A completed and better accessible account has been
given by Repka \cite{R}.\vspace{1.5mm}

For $l_j=-\frac12+i\lambda_j$\hspace{3cm}
$T_{l_1}\otimes T_{l_2}\sim
2\int\limits_{\R^+}^\oplus T_{-\frac12+i\lambda}d\lambda\,\oplus\,
\sum\limits_{l\in\Z}T_l$\,.
\\[.7mm]
For $l_1=-\frac12+i\lambda_1\,,\;l_2\in\N_0$\qquad\qquad\!\!
$T_{l_1}\otimes T_{l_2}\;\sim
\int\limits_{\R^+}^\oplus T_{-\frac12+i\lambda}d\lambda\oplus
\sum\limits_{l\ge0}T_l$\,.
\\[1.5mm]
For $l_j\in\N_0$\hspace{4.6cm}
$T_{l_1}\otimes T_{l_2}\ \sim \sum\limits_{l>l_1+l_2}T_l$\,.
\\
Similarly in the remaining cases.\vspace{1mm plus 1mm}

To get {\bf coefficients} for the unitary representations, we use
(corresponding to \cite{V}) an ortho{\it normal} basis $(e_m^l)$ of the
Hilbert space $\Cal H_l$ of $T_l$\,. For $l=-\frac12+i\lambda$
(principal series), we have $\Cal H_l=\Cal H$ and the basis is indexed by
$m\in\Z$\,. For $l\in\N_0$\,, the range is $m>l$ and for integers $l<0$\,:
$m\le l$\,.
\\
The basis vectors satisfy \;
$T_l(k_\varphi)\,e_m^l=e^{2mi\varphi}e_m^l=\chi_{2m}(k_\varphi)\,e_m^l$ \
("elliptic basis").
\\
We put \ $t_{mn}^l(g)=(T_l(g)e_n^l\mid e_m^l)$. This gives the
unitary matrix coefficients of $T_l(g)$.
\\
$t_{mn}^l$ is $(2m,2n)$-radial (we get only even integers, since we restrict
to representations of $\PSL(2,\R)$\,).
\\
For $l=-\frac12+i\lambda$\,, we have $t_{mn}^l\in B(G)$ for all $m,n\in\Z$
\,(it even belongs to the reduced Fourier-Stieltjes algebra $B_\rho(G)$,
i.e., the w*-closure of $A(G)$ in $B(G)$).
\\
For $l\in\Z$\,, the representations $T_l$ are {\it square-integrable},
thus $t_{mn}^l\in A(G)\cap L^2(G)$ for $l\in\N_0,\;m,n>l$ and for
$l<0,\;m,n\le l$\,.
\\
For $l=-\frac12+i\lambda$\,, the "non-radial component" of $t_{mn}^l$ is
described by \ $\fr P_{mn}^l(\ch2\tau)=t_{mn}^l
\begin{pmatrix}e^\tau & 0 \\  0 & e^{-\tau}
\end{pmatrix}$ for $\tau\ge0$ ($\ch$ denoting the hyperbolic cosine).
In \cite{V} the functions $\fr P_{mn}^l$ are defined (and investigated)
for all $l\in\C$\,, but (apart of the principal series)
using a non-normalized orthogonal basis for the matrix representation. For the
discrete series, the corresponding
functions arising from the {\it unitary} coefficients are denoted by
$\Cal P_{mn}^l$ in \cite{VK} ($l\in\Z$).
For $l\in\N_0\,,\;m,n>l$\vspace{1mm} they are related by
\linebreak
$\fr P_{mn}^l=
\Bigl(\dfrac{(m-l-1)!\,(n+l)!}{(m+l)!\,(n-l-1)!}\Bigr)^{\frac12}
\Cal P_{mn}^l$\,.

\medskip
Technically, the continuous part in the decomposition of tensor products
is more difficult to handle (and the appearance of multiplicities causes
additional complications). Therefore we restrict to the discrete
part.\vspace{1mm}

For $l_1=-\frac12+i\lambda\,,\;l_2\in\N_0$\,, we define the\vspace{.5mm}
{\it Clebsch-Gordan coefficients} by
\[
e_j^{l_1}\otimes e_m^{l_2}\,=
\,\sum\limits_{l\ge0}\,C(l_1,l_2,l;j,m,j+m)\,e_{j+m}^l\,
+\text{ cont.\,part .}
\]
The same for $l_1\in\Z$ with $l_1\ge-l_2-1$ \ (for $l_1<-l_2-1$ the discrete
part of $T_{l_1}\otimes T_{l_2}$ contains only $T_l$ with $l<0$\,).
We put $C(l_1,l_2,l;j,m,j+m)=0$ when $j+m\le l$ \ (in addition, for
$l_1\in\Z$\,, the coefficients will be $0$ outside the range $l>l_1+l_2$ for
$l_1\in\N_0$ and outside $0\le l\le l_1+l_2$ for $-1-l_2\le l_1<0$\,).
The isomorphism
between $T_l$ and a component of $T_{l_1}\otimes T_{l_2}$ is determined only
up to a factor of modulus $1$\,. This is fixed by requiring that \
$C(l_1,l_2,l;l-l_2,l_2+1,l+1)>0$ \ (of course, in the integer case this refers
only to those $l\ge0$ that have not been excluded above).
\\
For $l_1,l_2$ as above, this gives a decomposition of products in\vspace{-1mm}
$B(G)$
\begin{multline}\label{decpr}
t_{jj'}^{l_1}\;t_{mm'}^{l_2}=\\[.5mm]
\sum\limits_{l\ge0}
\overline{C(l_1,l_2,l;j,m,j+m)}\,C(l_1,l_2,l;j',m',j'+m')\,t_{j+m\,j'+m'}^l
+\text{ cont.\,part\,.} 
\end{multline}

Now, we consider the behaviour for {\bf large} $l_2$\,.
\begin{Pro}[Asymptotics of CG-coefficients]  \label{ascg}
For fixed $l_1=-\frac12+i\lambda\,,\;j,s\in\Z$ and finite $\kappa\ge 1$,
we have
\[
\lim\limits_{\substack{\phantom{a}l_2\to\infty \\[1.5pt]
\frac m{l_2}\,\to\kappa}}
C(l_1,l_2,l_2+s;j,m,j+m)=\fr P_{s\,j}^{l_1}(\kappa) \ .
\]
\end{Pro}
For $\kappa=1$, one has to add the restriction $m>l_2$\,.
Corresponding results hold for $l_1\in\Z$ (discrete series), e.g.,
when $l_1\in\N_0\,,\;j,s>l_1$\,, the limit is $\Cal P_{s\,j}^{l_1}(\kappa)$.
Similarly for the complementary series and unitary representations
of covering groups. This is the counterpart of a
classical result of Brussaard, Tolhoek (1957) on the CG-coefficients of
$SU(2)$.

Since $(\fr P_{s\,j}^{l_1}(\kappa))_{s,j\in\Z}$ is the matrix of a unitary
operator, its column vectors have norm $1$ (in $l^2(\Z)$). From
$\lVert e_j^{l_1}\otimes e_m^{l_2}\rVert=1$, it follows by orthogonality
that the norm of the continuous part in the decomposition of
$e_j^{l_1}\otimes e_m^{l_2}$ tends to $0$ for $l_2\to\infty$ (with\linebreak
$l_1,j$ fixed, $\frac m{l_2}\to\kappa$). The same holds for the
decomposition of $t_{jj'}^{l_1}\;t_{mm'}^{l_2}$ in \eqref{decpr}.
\\
It was already noted by Pukanszky that the densities arising in
the continuous part are given by analytic functions. Thus (with at most
contably many exceptions) all $\lambda\ge0$ will appear in the decomposition
of $e_j^{l_1}\otimes e_m^{l_2}$ (for $l_1=-\frac12+i\lambda_1$). But from
a more quantitative viewpoint, most of the product will be concentrated on
the (positive part of the) discrete series when $l_2$ is large.

\begin{proof}[Idea of Proof]
Recall the Fourier inversion formula:
\[
h(e)=\int\limits_0^\infty
\tr(T_{-\frac12+i\lambda}(h))\,\lambda\,\tah(\pi\lambda)\,d\lambda
+\sum\limits_{l\ge0}\,(l+\frac12)\bigl(\,\tr(T_l(h))+\tr(T_{-l-1}(h))\,\bigr)\;.
\]
for $h\in A(\PSL(2,\R))\cap L^1(\PSL(2,\R))$ and the extensions of the
representations to $L^1(\PSL(2,\R))$ for an appropriate choice of the
Haar measure. This describes also the Plancherel measure.
\\
On the level of coefficients, applied to $(2m,2n)$-radial functions
with $m,n\ge0$, this gives a generalization of the Mehler-Fock transformation
\[
g(x)=\sum\limits_{l=0}^{\min(m,n)-1}\,(l+\frac12)\,b(l)\,\Cal P_{mn}^l(x)
+\int\limits_0^\infty\,a(\lambda)\,\fr P_{mn}^{-\frac12+i\lambda}(x)
\,\lambda\,\tah(\pi\lambda)\,d\lambda
\]
with $b(l)=\int\limits_1^\infty\,g(x)\,\Cal P_{mn}^l(x)\,dx$ \ for
$g\in L^2([1,\infty])$ \ (convergence in $L^2$). Thus the discrete part is
just the expansion with respect to the orthogonal system
$(\Cal P_{mn}^l)\subseteq L^2([1,\infty])$ ($m,n$ fixed) and the coefficients
are obtained from inner products.
\\
We apply this to \,$g=\fr P_{ss}^{l_1}\,\Cal P_{l_2+1\,l_2+1}^{l_2}$ and get
for $l=l_2+s$ by \eqref{decpr}
\begin{multline*}
\lvert\, C(l_1,l_2,l_2+s;s,l_2+1,l_2+s+1)\rvert^2=\\
(l_2+s+\frac12)\,\int\limits_1^\infty\,\fr P_{ss}^{l_1}(x)\,
\Cal P_{l_2+1\,l_2+1}^{l_2}(x)\,\Cal P_{l_2+s+1\,l_2+s+1}^{l_2+s}(x)\,dx
\end{multline*}
By \cite{V} we have
$\Cal P_{l+1\,l+1}^l(x)=\fr P_{l+1\,l+1}^l(x)=\Bigl(\dfrac2{x+1}\Bigr)^{l+1}$.
\vspace{1mm}It follows easily that for $l_2\to\infty$ and $s\in\Z$ fixed,
$(l_2+s+\frac12)\,\Cal P_{l_2+1\,l_2+1}^{l_2}\,
\Cal P_{l_2+s+1\,l_2+s+1}^{l_2+s}\to\delta_1$ (point measure)
holds weakly with respect to
bounded continuous functions on $[1,\infty[$\,. Since $\fr P_{ss}^{l_1}(1)=1$,
this gives $\lvert\, C(l_1,l_2,l_2+s;s,l_2+1,l_2+s+1)\rvert\to1$ (when
$l_1=-\frac12+i\lambda$ is fixed) and by our choice of the phase, we get
$C(l_1,l_2,l_2+s;s,l_2+1,l_2+s+1)\to1$.
\\[1mm]
Next we take $g=\fr P_{s\,j}^{l_1}\,\Cal P_{l_2+1\,m}^{l_2}$ and get
for $l=l_2+s$ by \eqref{decpr}
\begin{multline*}
\underset{\textstyle\to1^{\phantom{l}}}
{\overline{C(l_1,l_2,l_2+s;s,l_2+1,l_2+s+1)}}\:
C(l_1,l_2,l_2+s;j,m,j+m)=\\[-3mm]
(l_2+s+\frac12)\,\int\limits_1^\infty\,\fr P_{s\,j}^{l_1}(x)\,
\Cal P_{l_2+1\,m}^{l_2}(x)\,\Cal P_{l_2+s+1\,j+m}^{l_2+s}(x)\,dx
\end{multline*}
Let $\mu_{l_2\mspace{.5mu}m}$ be the measure on $[1,\infty[$ with density
$(l_2+s+\frac12)\,\Cal P_{l_2+1\,m}^{l_2}\,\Cal P_{l_2+s+1\,j+m}^{l_2+s}$.
Again one can use the formulas of \cite{V} for $\fr P_{l+1\,m}^l(x)$.
With a slight change of coordinates, one gets that
$\dfrac{\mu_{l_2\mspace{.5mu}m}}{\lVert\mu_{l_2\mspace{.5mu}m}\rVert}$ has a
$\beta'$-distribution and from the values of expectation and variance one can
conclude that $\lVert\mu_{l_2\mspace{.5mu}m}\rVert\to1$ and
$\mu_{l_2\mspace{.5mu}m}\to\delta_\kappa$ for
$l_2\to\infty,\;\dfrac m{l_2}\to\kappa$\,.
\end{proof}

In the next step we use {\bf ultraproducts} to work with these limit
relations. 
Such constructions for group representations have been done by Cowling and
Fendler.
\\[1mm]
We take some element $p\in\beta\N\setminus\N$
\,(Stone-\v Cech compactification).
The ultraproduct of the Hilbert spaces $(\Cal H_l)_{l>0}$
\,(with respect to $p$)
is denoted by $\Cal H_p$\,. It consists of equivalence classes of all
sequences $(h_l)\in\prod\Cal H_l$ such that
$\lim\limits_{l\to p}\,\lVert h_l\rVert<\infty$\,, factoring by the subspace
of sequences with $\lim\limits_{l\to p}\,\lVert h_l\rVert=0$\,.
We use the notation \ $\lim\limits_{l\to p}\,h_l$ \;to denote the
equivalence class of $(h_l)$. $\Cal H_p$ is again a Hilbert space and we
get a representation $T_p$ of the C*-algebra $\VN(G)$ on $\Cal H_p$ putting
\;$T_p(S)(\lim\limits_{l\to p}\,h_l)=\lim\limits_{l\to p}\,T_l(S)h_l$
(for $S\in \VN(G)$\,).
\\
Each function $f\!:\N\to\N$ satisfying $f(l)>l\ \;\forall\,l$
(or more generally, \,$\lim\limits_{l\to p}\,f(l)-l>0$\,) defines a
unit vector in $\Cal H_p$ by \,$e(p,f)=\lim\limits_{l\to p}\,e_{f(l)}^l$\,.
Of course, it is enough to require that $f$ is defined
for $l\ge l_0$\,. For functions $f,f'$ we get a coefficient functional by
\;$(t_{ff'}^p\,,\,S\,)=\bigl(T_p(S)\,e(p,f')\mid e(p,f)\bigr)$
for $S\in \VN(G)$\,. Then $t_{ff'}^p\in \VN(G)'$ (dual space) and
$t_{ff'}^p=\lim\limits_{l\to p}\,t_{f(l)\,f'(l)}^l$ (w*-limit).
\\
Recall that $\beta\N\setminus\N$ is a $\Z$-module under addition.
Thus we get in the same way Hilbert spaces $\Cal H_{p+s}$ and
representations $T_{p+s}$ for all $s\in\Z$\,.
\\
For $f$ as above, put \,$\kappa_p(f)=\lim\limits_{l\to p}\,\dfrac{f(l)}l$
(possibly infinite).
\\
Write $\kappa=\kappa_p(f),\;\kappa'=\kappa_p(f')$. Assuming,
$1<\kappa,\kappa'<\infty,\ l_1=-\frac12+i\lambda$\,, we get from
\eqref{decpr} and Proposition\;\ref{ascg} \
\[
t_{jj'}^{l_1}\odot t_{ff'}^p=
\lim\limits_{l_2\to p}\,t_{jj'}^{l_1}\;t_{f(l_2)f'(l_2)}^{l_2}=
\sum\limits_{s\in\Z}\,
\overline{\fr P_{sj}^{l_1}(\kappa)}\,\fr P_{sj'}^{l_1}(\kappa')\,
\lim\limits_{l_2\to p}\,t_{f(l_2)+j\:f'(l_2)+j'}^{l_2+s}
\]
(note that
$\bigl(\,\overline{\fr P_{sj}^{l_1}(\kappa)}\,
\fr P_{sj'}^{l_1}(\kappa')\,\bigr)_{s\in\Z}
\in l^1$\,). Put $u(l)=l-1$ for $l\in\Z$\,, then\linebreak
$\lim\limits_{l_2\to p}\,t_{f(l_2)+j\:f'(l_2)+j'}^{l_2+s}=
t_{f\circ u^s+j\:f'\circ u^s+j'}^{p+s}$ and we arrive at
\begin{equation}\label{declim}
t_{jj'}^{l_1}\odot t_{ff'}^p=
\sum\limits_{s\in\Z}\,
\overline{\fr P_{sj}^{l_1}(\kappa)}\,\fr P_{sj'}^{l_1}(\kappa')\,
t_{f\circ u^s+j\:f'\circ u^s+j'}^{p+s}\ .
\end{equation}
Next, we consider \,
$\overline{\Cal H}_p=\bigoplus\limits_{s\in\Z}\,\Cal H_{p+s}$ \
($l^2$-sum) and the corresponding representation
$\overline T_p=\bigoplus\limits_{s\in\Z}\,T_{p+s}$ \,of $\VN(G)$.
\\
For $1<\kappa<\infty$\,, $\Cal K_\kappa$ shall be the closed subspace of
$\Cal H_p$ generated by the vectors $e(p,f)$, taking all functions $f$ with
$\kappa_p(f)=\kappa$\,. We put
\,$\Cal K=\bigoplus\limits_{1<\kappa<\infty}\Cal K_\kappa$\,.
\\
$U(\lim\limits_{l\to p+s}\,h_l)=\lim\limits_{l\to p+s+1}\,h_{l-1}$
\,defines an isometric isomorphism of $\Cal H_{p+s}$ and $\Cal H_{p+s+1}$
and this extends to a unitary operator
$U\!:\overline{\Cal H}_p\to\overline{\Cal H}_p$ \;(in particular
$U\bigl(e(p+s,f)\bigr)=e(p+s+1,f\circ u)$\,).
Let $\overline{\Cal K}_\kappa$ be the closed $U$-invariant subspace of
$\overline{\Cal H}_p$ generated by $\Cal K_\kappa$ \,(it is generated by
the vectors $e(p+s,f)$, taking all functions $f$ with
$\kappa_{p+s}(f)=\kappa$ for some $s\in\Z$\,). Clearly,
$\overline{\Cal K}_\kappa\perp\overline{\Cal K}_{\kappa'}$ holds for
$\kappa\ne\kappa'$ and we write
$\overline{\Cal K}=\bigoplus\limits_{1<\kappa<\infty}\overline{\Cal K}_\kappa$
\,(the closed $U$-invariant subspace of
$\overline{\Cal H}_p$ generated by $\Cal K$\,).
\;$V\bigl(e(p+s,f)\bigr)=e(p+s,f+1)$ defines a unitary operator on
$\overline{\Cal K}_\kappa$ (for $1<\kappa<\infty$) and this extends to a
unitary operator \,$V\!:\overline{\Cal K}\to\overline{\Cal K}$ \,satisfying
\ $V(\Cal K_\kappa)\subseteq\Cal K_\kappa$ and
$VU=UV$ on $\overline{\Cal K}$\,. \ (For $\kappa=1, \ V$ is no longer
surjective).
\\[2mm]
For a fixed function $f$ with $\kappa=\kappa_p(f)$ satisfying
$1<\kappa<\infty$\,, it follows easily that \
$\{e(p+s,f\circ u^s+j)\}=\{U^sV^je(p,f): s,j\in\Z\}$ \;defines an
orthonormal system of vectors in $\overline{\Cal K}_\kappa$\,.
\\
A special case, used below, will be the functions $f_\kappa(l)=[\kappa\,l]$
\,(integer part), satisfying $\kappa_p(f_\kappa)=\kappa$ \,for each $p$
and $1<\kappa<\infty$\,.\vspace{-.7mm}

\begin{Lem}
For $\lambda\in\R\,,\;j\in\Z\,,\;1<\kappa<\infty$
\\
$A_j^\lambda=
V^j\,\sum\limits_{s\in\Z}\fr P_{sj}^{-\frac12+i\lambda}(\kappa)\:
\bigl\lvert 2s\bigr\rvert^{i\lambda}\,U^s$ \ defines a bounded linear
operator \
$\Cal K_\kappa\to\overline{\Cal K}_\kappa$\,.
\\
Taking $A_j^\lambda=0$ on $\Cal K^\perp$ (in particular, $A_j^\lambda=0$ on
$\Cal H_{p+s}$ when $s\neq0$) gives a bounded linear operator
$A_j^\lambda\!:\overline{\Cal H}_p\to\overline{\Cal H}_p$ \,satisfying
\,$V\mspace{-1mu}A_j^\lambda=A_j^\lambda\mspace{1mu}V$ on $\overline{\Cal K}$.
\end{Lem}
(Here we adopt $0^{i\lambda}=1$).

\begin{Cor}
Given $e,e'\in\Cal K$ define $t\in \VN(G)'$ by
\,$(t,S)=(T_p(S)\,e'\mid e\,)$. Then for $l=-\frac12+i\lambda \ (\lambda\in\R)$
and $j,j'\in\Z$ we have \ $(t_{jj'}^l\odot t,S)=
(\overline T_p(S)A_{j'}^\lambda e'\mid A_j^\lambda e)$ \ ($S\in \VN(G)$).
\end{Cor}

\begin{Lem}  \label{dense}
$\overline T_p(\VN(G))$ is w*-dense in
\;$\prod\limits_{s\in\Z}\Cal B(\Cal H_{p+s})$\,.
\end{Lem}

In particular, this implies that $T_p$ is irreducible and $(T_p,\Cal H_p)$
is the cyclic representation for the state $t_{ff}^p$
(with cyclic vector $e(p,f)$\,) for every function $f$ as above.
Furthermore (slightly more general than in Lemma\;\ref{dense}), one has
\,$T_p\nsim T_{p'}$ for $p\ne p'$. Considering $L^1(G)$ as a (w*-dense)
subalgebra of $\VN(G)$, it is not hard to see that $T_p(h)=0$ for
$h\in L^1(G)$, hence these are singular representations of $\VN(G)$.
\\[2mm]
For the final step we need a refinement of Lemma\;\ref{dense}.
Although $\overline T_p(\VN(G))$ is not a von Neumann algebra, the fact that
$\VN(G)$ is a von Neumann algebra allows to get a stronger result on the
size of \,$\overline T_p(\VN(G))$.
\\
Recall that the representations $T_l$ are {\it square integrable} for
$l\in\Z$\,. Thus they are equivalent to subrepresentations of the left regular
representation on $L^2(G)$ and we can consider
$\prod\limits_{l\ge0}\Cal B(\Cal H_l)$ as a subalgebra of $\VN(G)$.
\\
For $1\le\alpha<\beta\le\infty$ let $P_{\alpha\beta}\in \VN(G)$ be the
orthogonal projection on the closed subspace of
$\bigoplus\limits_{l>0}\Cal H_l$ generated by
\,$\bigl\{e_m^l: \alpha<\dfrac ml<\beta,\;l>0\bigr\}$. For
\,$\alpha<\beta\le \alpha'<\beta'$, \,it follows that
$P_{\alpha\beta}P_{\alpha'\beta'}=P_{\alpha'\beta'}P_{\alpha\beta}=0$\,.
For $\alpha<\kappa<\beta$ we have
\,$\overline{\Cal K}_\kappa\subseteq
\im\bigl(\overline T_p(P_{\alpha\beta})\bigr)$\,.

\begin{Lem}  \label{dense2}
Assume that $\alpha_m\nearrow\infty$\,. For $m\ge1$,
\begin{itemize}\itemindent-5mm
\item[]\,$E_m\ (\subseteq\overline{\Cal H}_p)$ shall be a finite dimensional
subspace of \;$\im\bigl(\overline T_p(P_{\alpha_m\alpha_{m+1}})\bigr)$\,,
\item[]\,$S_m\in\Cal B(\overline{\Cal H}_p)$ are such that
\,$\lVert S_m\rVert\le1$,
\;$S_m(E_m)\subseteq \im\bigl(\overline T_p(P_{\alpha_m\alpha_{m+1}})\bigr)$
\ and\linebreak
$S_m(\Cal H_{p+s})\subseteq\Cal H_{p+s}$ \,for all $s\in\Z$\,.
\end{itemize}
Then there exists \,$S\in \VN(G)$ such that
\,$\bigl\Vert\,\bigl(S_m-\overline T_p(S)\bigr){\big\vert} E_m\,\bigr\Vert\to0$ for
$m\to\infty$\,.
\end{Lem}
\medskip

At the Harmonic Analysis Conference in Istanbul 2004, I talked about
the case $G=SU(2)$. For that group, one could use a limit of averages of
states $t_{ff}^p$ (for\linebreak
$f=f_\kappa$\,; approaching Lebesgue measure on $[-1,1]$\,)
to get a singular state\linebreak $\zeta\in \VN(G)'$
satisfying $\lVert f\odot\zeta\rVert=\lVert f\rVert$ for all $f\in A(G)$.
This cannot exist for $G=\SL(2,\R)$, because of non-amenability.
Instead of this, we will use another type of asymptotics.

Now, we fix $p\in\beta\N\setminus\N$ and write $\overline T$  for
$\overline T_p$\,. We choose
$p_1\in\beta\N\setminus\N$
satisfying $(2^m)\in p_1$  (a sufficiently "thin" ultrafilter).
$(\overline{\Cal H}_p)_{p_1}$ shall denote the ultrapower of
$\overline{\Cal H}_p$ with respect to $p_1$\,. If $(h^{(n)})$ is a
bounded sequence in $\overline{\Cal H}_p$\,, we write, as before,
$\lim\limits_{n\to p_1}h^{(n)}$ for the corresponding equivalence class,
defining an element of $(\overline{\Cal H}_p)_{p_1}$\,. The representation
$\overline T$ of $\VN(G)$ on $\overline{\Cal H}_p$ defines a representation
$\overline{\overline T}$ of $\VN(G)$ on $(\overline{\Cal H}_p)_{p_1}$\,.
We define $\bar{\bar e}\in(\Cal K)_{p_1}\subseteq
(\Cal H_p)_{p_1}$ and $\zeta\in \VN(G)'$
by
\[
\boxed{\quad\bar{\bar e}=
\lim\limits_{n\to p_1}\frac1n\sum\limits_{r=1}^{n^2-1}
e(p,f_{\ch(n+{\textstyle\frac rn})})\ ,\qquad
(\zeta,S)=(\,\overline{\overline T}(S)\,\bar{\bar e}\mid \bar{\bar e}\,)
\quad}\vspace{1mm}
\]
$V$ defines a unitary
operator on $(\overline{\Cal K})_{p_1}$\,, again denoted by $V$\,.
Since (for fixed $f$)
$\{V^je(p,f):j\in\Z\}$ is an orthonormal family in $\Cal K_{\kappa(f)}$,
\,$\Cal K_\kappa\perp\Cal K_{\kappa'}$ for $\kappa\neq\kappa'$
and $V$ is unitary\,, it follows that
$\{V^j\bar{\bar e}:j\in\Z\}$ is orthonormal in
$(\Cal K)_{p_1}$\,. Thus for \,$\gamma=(\gamma_n)\in l^2(\Z)$,
\;$\bar{\bar e}(\gamma)=\sum_{j\in\Z}\gamma_jV^j\bar{\bar e}$ defines an
isometric embedding $l^2(\Z)\to(\Cal K)_{p_1}$\,.
For $\gamma_1,\gamma_2\in l^2(\Z)$, we define
$(\zeta_{\gamma_1\gamma_2}\,,S)=
(\,\overline{\overline T}(S)\,\bar{\bar e}(\gamma_2)\mid
\bar{\bar e}(\gamma_1)\,)$ and for brevity\linebreak (\,$(e_n)$
denoting the standard basis of $l^2(\Z)$) \,$\zeta_{ij}=\zeta_{e_ie_j}$
\,(thus $\bar{\bar e}=\bar{\bar e}(e_0),\ \zeta=\zeta_{00}$).
\\
For $g\in\Cal K(\,(\R\setminus\{0\})\times\Z\,)$ \
($\Cal K(\Omega)$\,: continuous
functions with compact support), we put\vspace{-4mm}
\[
\varphi(g)=
\lim\limits_{n\to p_1}\frac1n\sum\limits_{r=1}^{n^2-1}
\sum\limits_{j,s\in\Z}\,g\bigl(\frac{2s}{e^c},j\bigr)\,(-1)^s
\frac{\sqrt 2}{e^{c/2}}\ U^sV^j\,e(p,f_{\ch c})\qquad
\text{with}\quad c=n+\frac rn
\]
Note that the support condition makes the sum finite and restricts it to
$s\ne0$, hence $\varphi(g)\perp(\Cal H_p)_{p_1}$\,.

\begin{Lem}  \label{isom}
$\varphi(g)\in(\overline{\Cal K})_{p_1}\subseteq
(\overline{\Cal H}_p)_{p_1}$\,, \
$\lVert\varphi(g)\rVert=\lVert g\rVert_2$\,.
\\
Thus $\varphi$ extends to an isometry
\;$\varphi\!:L^2(\R\times\Z)\to(\overline{\Cal H}_p)_{p_1}$\,.
\\
Putting $\varphi_1(g+\gamma)=\varphi(g)+\bar{\bar e}(\gamma)$ defines an
isometry 
\,$\varphi_1\!:L^2(\R\times\Z)\oplus l^2(\Z)\to(\overline{\Cal H}_p)_{p_1}$\,.
\end{Lem}

Let $P\in\Cal B\bigl((\overline{\Cal H}_p)_{p_1}\bigr)$ be the
orthogonal projection to \,$\varphi(L^2(\R\times\Z))$\,.
For \linebreak$S\in \VN(G),\;g,h\in L^2(\R\times\Z)$ put
\,$(\psi(S)g\mid h)=(\overline{\overline T}(S)\varphi(g)\mid\varphi(h))$\,.
This defines a contractive linear mapping
$\psi\!:\VN(G)\to\Cal B(L^2(\R\times\Z))$\,, \ $\psi(\VN(G))$ being
isometrically isomorphic to the dilation
$P\,\overline{\overline T}(\VN(G))\,P$\,.
\\
Similarly, for $P_1$ the projection to
$\varphi_1(L^2(\R\times\Z)\oplus l^2(\Z))$, one
gets \,$\psi_1\!:\VN(G)\to\linebreak
\Cal B(L^2(\R\times\Z))\oplus\Cal B(l^2(\Z))$ \
(note that $(\Cal H_p)_{p_1}$ is invariant under
$\overline{\overline T}(\VN(G))$\,).
\\[2mm]
For $n=2^m,\;\alpha_m=\ch2^m$, the $n$-th term in the limits defining 
$\bar{\bar e}$ and $\varphi(g)$ belongs to
$\im\bigl(\overline T_p(P_{\alpha_m\alpha_{m+1}})\bigr)$. This makes it
possible to apply Lemma \ref{dense2}.\vspace{-2mm}

\begin{Lem}  \label{dense3}
$\psi(\VN(G))$ is w*-dense in \,$\Cal B\bigl(L^2(]-\infty,0]\times\Z)\bigr)
\oplus\Cal B\bigl(L^2([0,\infty[\times\Z)\bigr)$.
\end{Lem}
Similarly, for $\psi_1$ one has to add a sum with $\Cal B(l^2(\Z))$. As above, the
w*-closure of $\psi(\VN(G))$ is isometrically isomorphic to
$P\,\overline{\overline T}{(\VN(G))}^-P$ \ (${}^-$ denoting the\linebreak
w*-closure
in $\Cal B\bigl((\overline{\Cal H}_p)_{p_1}\bigr)$\,). The proof shows
that corresponding density results hold for the image
of the unit ball of $\VN(G)$.
\medskip

For the final step, we will use the {\it Whittaker functions}. They
are defined by
\[
W_{\lambda,\mu}(z)=
\frac{z^{\mu+\frac12}\,e^{-\frac z2}}{\Gamma(\mu-\lambda+\frac12)}\,
\int\limits_0^\infty e^{-zu}\,u ^{\mu-\lambda-\frac12}\,
(1+u)^{\mu+\lambda-\frac12}\,du
\]
for \,$\re z>0$, $\re(\mu-\lambda+\frac12)>0$ and then for all
$\lambda,\mu\in\C$ by analytic continuation.

\begin{Pro}[Approximation of coefficients]  \label{apco}
For \,$n\in\Z$, \,$l=-\frac12+i\lambda$ fixed,
\[
\sup\,\biggl\{\,\Bigl\lvert\,\fr P_{mn}^l(\ch \tau)-
\frac{(-1)^{n-m}}{m^{l+1}\Gamma(n-l)}\,
W_{n,i\lambda}\Bigl(\frac{4m}{e^\tau}\Bigr)\,\Bigr\rvert\:
e^{\textstyle\frac\tau2}\,m^2:\;\tau\ge0,\,m\ge n,\,m>0\biggr\}
\]
is finite.
\end{Pro}\vspace{-2mm}
In particular,\quad
$
\lim\limits_{m\to\infty}\,\Bigl(\,\fr P_{mn}^l(\ch \tau)-
\dfrac{(-1)^{n-m}}{m^{l+1}\Gamma(n-l)}\,
W_{n,i\lambda}\bigl(\frac{4m}{e^\tau}\bigr)\,\Bigr)\,
e^{\textstyle\frac\tau2}=0
$\vspace{1mm}
holds uniformly for $\tau\ge0$\,. \
This complements classical results on the asymptotic behaviour of
$\fr P_{mn}^l$ for fixed values
$l,m,n$\,; \ e.g., if $m=n,\ \lambda\ne0$ \,recall that\linebreak
\,$\fr P_{mm}^l(\ch\tau)\,e^{\textstyle\frac\tau2}-
\dfrac2{\sqrt{\pi\lambda\tah(\pi\lambda)}}\cos(\lambda\tau+\eta)\to0$
\ for $\tau\to\infty$ (where $\eta\in\R$ depends on $\lambda$ and $m$).\\
Moreover, the Proposition implies also that the row vector
$\bigl(\,\fr P_{mn}^l(\ch \tau)\bigl)_{m>0}$ can be approximated
in $l^2$-norm by $\Bigl(\frac{(-1)^{n-m}}{m^{l+1}\Gamma(n-l)}\,
W_{n,i\lambda}\bigl(\frac{4m}{e^\tau}\bigr)\Bigr)$ for $\tau\to\infty$\,.
An approximation for the "lower half"
$\bigl(\,\fr P_{mn}^l(\ch \tau)\bigl)_{m<0}$ is obtained using the
identity $\fr P_{mn}^l=\fr P_{-m\,-n}^l$\,.\vspace{1mm}

For $j\in\Z,\;\lambda\in\R,\;l=-\frac12+i\lambda$\,, we put
\[
g_{j,\lambda}(x,j')=\begin{cases}
\hspace{1.6cm} 0&\text{for }j'\ne j\vspace{.5mm} \\
\ \dfrac{(-1)^j\,2^{i\lambda}}{\Gamma(j-l)\,\sqrt{\;x}}\
W_{j,i\lambda}(2x)&\text{for } j'=j\,,\;x>0\vspace{1.5mm} \\
\dfrac{(-1)^j\,2^{i\lambda}}{\Gamma(-j-l)\sqrt{-x}}\,
W_{-j,i\lambda}(-2x)&\text{for } j'=j\,,\;x<0
\end{cases}\]\vspace{-3mm}
Then $g_{j,\lambda}\in L^2(\R\times\Z)$.
\\[5mm]
$A_j^\lambda\in\Cal B(\overline{\Cal H}_p)$ defines a bounded
operator on $(\overline{\Cal H}_p)_{p_1}$\,, again denoted by
$A_j^\lambda$\,. At the other side,
for $g\in L^2(\R\times\Z)$, we define $(Vg)(t,j)=g(t,j-1)$.\vspace{-2mm}
\begin{Lem}  \label{ebar}
We have \ $A_j^\lambda\, \bar{\bar e}=\varphi(g_{j,\lambda})$\,.
\\
Furthermore, $V\varphi(g)=\varphi(Vg)$ for $g\in L^2(\R\times\Z)$,
in particular, $A_j^\lambda V^s\, \bar{\bar e}=\varphi(V^sg_{j,\lambda})$
holds for all $s\in\Z$\,.
\end{Lem}

\begin{Cor}
$(t_{jj'}^l\odot\zeta\,,S)=
(\psi(S)\,g_{j',\lambda}\mid g_{j,\lambda})$ \ ($S\in \VN(G)$).
\\[1mm]
More generally, \;$(t_{jj'}^l\odot\zeta_{ii'}\,,S)=
(\psi(S)\,V^{i'}g_{j',\lambda}\mid V^ig_{j,\lambda})$
\;for $i,i'\in\Z$\,.\vspace{2mm}
\end{Cor}

The basis of $L^2(\R)$ used by \cite{V} to define the coefficients of
$T_l$ for $l=-\frac12+i\lambda$ is given by
\,$e_m^l(x)=\dfrac{(-1)^m}{\sqrt{\pi}}\,e^{2mi\arctan(x)}(1+x^2)^l=
\dfrac1{\sqrt{\pi}}\Bigl(\dfrac{x-i}{x+i}\Bigr)^m\bigl(1+x^2\bigr)^l$
.\vspace{1mm}

We consider the real Fourier transform
$\hat f(y)=\dfrac1{\sqrt{2\pi}}\,\int\limits_{\R}e^{-ixy}f(x)\,dx$\,.
\vspace{-2mm}Then we have
\[
\widehat{e_m^l}(y)=(-1)^m\,
\frac{2^{i\lambda}\lvert y\rvert^{-\frac12-i\lambda}}{\Gamma(\sgn(y)m-l)}
\,W_{\sgn(y)m,i\lambda}(2\lvert y\rvert)=
g_{m,\lambda}(y,m)\,\lvert y\rvert^{-i\lambda}\ .
\]
(The functions $e_m^l$ are not integrable, so strictly speaking, this is the
Fourier-Plancherel transform).
\\
For $h=\begin{pmatrix}a & 0\\b& \frac1a
\end{pmatrix}\in H$\,, we have
$(T_l(h)f)(x)=\lvert a\rvert^{-2l}f(a^2x+a\,b)$.
Composition with Fourier transform defines equivalent representations
(Whittaker model) $\pi_\lambda(g)\hat f=(T_l(g)f)\sphat$\;.
For $h\in H$ this gives \,$(\pi_\lambda(h)\,\eta)(y)=
\lvert a\rvert^{-1-2i\lambda}
\,e^{iy{\textstyle\frac ba}}\,\eta\bigl(\dfrac y{a^2}\bigr)$.
\\[.5mm]
Put \,$(\rho_\lambda\eta)(y)=\lvert y\rvert^{i\lambda}\,\eta(y)$\,.
Then $\rho_\lambda\!:L^2(\R)\to L^2(\R)$ is an isometric isomorphism
and \,$\pi_\lambda(h)=\rho_\lambda^{-1}\circ\pi_0(h)\circ\rho_\lambda$
\,(in particular, all $T_l$ and $\pi_\lambda$ define equivalent
representations of $H$\,).
$\pi_0$ splits into two irreducible representations (the
restrictions to $L^2(]-\infty,0]$ and $L^2([0,\infty[$) and these are
the only infinite dimensional irreducible unitary representations of $H$
\,(up to equivalence). Thus $\pi_0$ defines a normal isomorphism of the
von Neumann algebras $\VN(H)$ and $\Cal B\bigl(L^2(]-\infty,0])\bigr)\oplus
\Cal B\bigl(L^2([0,\infty[)\bigr)$ \,and this extends to a normal isomorphism
$\tilde\pi_0$ of the
von Neumann algebras $\Cal A=\VN(H)\bar\otimes \Cal B(l^2(2\Z))$ and
$\Cal B\bigl(L^2(]-\infty,0]\times\Z)\bigr)\oplus
\Cal B\bigl(L^2([0,\infty[\times\Z)\bigr)$.
\\
We have $g_{j,\lambda}(\cdot,j)=\rho_\lambda\,\widehat{e_j^l}$\,,
consequently $\pi_0(S)\,g_{j,\lambda}(\cdot,j)=
\rho_\lambda\bigl(\pi_\lambda(S)\,\widehat{e_j^l}\bigr)=\linebreak
\rho_\lambda\bigl((T_l(S)\,e_j^l)\sphat\,\bigr)$, resulting in
\[\tag{3}\label{pi0}
\bigl(\pi_0(S)\,g_{j',\lambda}(\cdot,j')\mid g_{j,\lambda}(\cdot,j)\bigr)=
(\,S\,,t_{jj'}^l\vert H) \text{\quad for\quad} S\in \VN(H)\ .\vspace{1mm}
\]
For $f\in M(A(G))$ put $\Phi(f)=(f_{mn})_{m,n\in2\Z}$ with
$f_{mn}=(\chi_n*f*\chi_m)\vert H$ \,(the matrix operators, used in the
Theorem). 
Put
$\lambda=\lim_{x\to\infty}f_{00}(x),\ f_0=f-\lambda,\
\Phi_1(f)=\Phi(f_0)+\lambda\,e_{00}$\,. Extend $\tilde\pi_0$ to a normal
isomorphism $(\pi_0\oplus1)\sptilde$ of the
von Neumann algebras $(\VN(H)\oplus\C)\bar\otimes \Cal B(l^2(2\Z))$ and
$\Cal B\bigl(L^2(]-\infty,0]\times\Z)\bigr)\oplus
\Cal B\bigl(L^2([0,\infty[\times\Z)\bigr)\oplus\Cal B(l^2(\Z))$.
Recall that $f$ is said to be $K$-finite, if only finitely many $f_{mn}$ are
non-zero.

\begin{Lem}  \label{pred}
For $f\in M(A(G))\cap C_0(G),\ \;\Phi(f)$ defines an element of the predual
of $\VN(H)\bar\otimes \Cal B(l^2(2\Z))$ and, if $f$ is $K$-finite, we have
\begin{gather*}\tag{4}\label{theta}
(f\odot\zeta\,,S)=
\bigl(\,\tilde\pi_0^{-1}\circ\psi(S)\,,\,\Phi(f)\,\bigr) \
\text{ for }S\in \VN(G)\ ,
\\
(f\odot\zeta_{i'i}\,,S)=
\bigl(\,\tilde\pi_0^{-1}\circ\psi(S)\,,\,\theta_f(e_{2i\,2i'})\,\bigr) \
\text{ for }i,i'\in\Z\ .
\end{gather*}
$\theta_f$ extends to a continuous linear mapping
$\Cal N(l^2(2\Z))\to\Cal A_*$ with
$\lVert \theta_f\rVert\le\lVert f\rVert_M$\,.
\\[1mm]
For general $f\in M(A(G)), \;f_0\in C_0(G)$ holds and  $\Phi_1(f)$ defines an
element of the predual of
\,$\bigl(\VN(H)\oplus\C\bigr)\bar\otimes \Cal B(l^2(2\Z))$\,.
If $f$ is $K$-finite, we have
\[
(f\odot\zeta\,,S)=
\bigl(\,((\pi_0\oplus1)\sptilde)^{-1}\circ\psi_1(S)\,,\,\Phi_1(f)\,\bigr) \
\text{ for }S\in \VN(G)\ .
\]
With $\theta_{1f}=\theta_{f_0}+\lambda$
\;(\,$\theta_{1f}\!:\Cal N(l^2(2\Z))\to\Cal A_*\oplus\Cal N(l^2(2\Z))$,
identified with the predual of \,$\Cal A\oplus \Cal B(l^2(2\Z))$), we get
(for $K$-finite $f$)
\[
(f\odot\zeta_{i'i}\,,S)=
\bigl(\,((\pi_0\oplus1)\sptilde)^{-1}\circ\psi_1(S)\,,
\,\theta_{1f}(e_{2i\,2i'})\,\bigr) \
\text{ for }S\in \VN(G)\ .
\]
\end{Lem}
\begin{Cor}
$\lVert\Phi_1(f)\rVert=\lVert\Phi(f_0)\rVert+\lvert\lambda\rvert\le
\lVert f\odot\zeta\rVert$ \ and \
$\lVert \theta_{f_0}\rVert+\lvert\lambda\rvert=
\lVert \theta_{1f}\rVert\le\lVert f\rVert_M$
\;holds for all $f\in M(A(G))$.
\end{Cor}
As indicated earlier this supplies the remaining step for the proof of the
Theorem. 

\begin{proof}[Idea of Proof]
Recall that the left and right actions of $G$ on $A(G)$ are continuous
and isometric. It follows easily that $f\in M(A(G))$ implies
$\mu*f,\,f*\mu\in M(A(G))$ for every bounded measure $\mu$ \,on $G$\,.
\\
For general $f\in M(A(G))$, the same argument as in \cite{CH} gives
$f\vert H\in B(H)$.
As mentioned before, the unitary dual of $H$ ($ax+b$\,-group) has a very
simple structure and this implies $B(H)=A(H)+B(H/[H,H])$. Thus for
$f\in M(A(G))\cap C_0(G)$,
we get (since $[H,H]$ is not compact) $f\vert H\in A(H)$, in particular,
$f_{mn}\in A(H)$ for all $m,n\in\Z$..

Let $M_1$ be the set of all $f\in M(A(G))\cap C_0(G)$ such that \eqref{theta}
holds. For $\gamma_1,\gamma_2\in l^2(\Z)$, it follows from the definition
that
$\lVert\zeta_{\gamma_1\gamma_2}\rVert\le
\lVert\gamma_1\rVert\,\lVert\gamma_2\rVert$ and this gives\linebreak
$\lVert\sum_{k,l}\alpha_{kl}\zeta_{kl}\rVert\le
\lVert(\alpha_{kl})\rVert_{\Cal N}$ for $(\alpha_{kl})\in\Cal N(l^2(\Z))$.
Thus if $f\in M_1$\,, then \eqref{theta} implies, using bilinearity
that $\theta_f$ is bounded and $\lVert\theta_f\rVert\le\lVert f\rVert_M$
and then the earlier argument, based on Proposition 1 shows
$\lVert f\rVert_M=\lVert f\rVert_{M_0}=\lVert \theta_f\rVert$. If
$f\in M_1$ is $(m,n)$-radial, it follows that
$\lVert f\rVert_M=\lVert \Phi_f\rVert=\lVert f\vert H\rVert_{A(H)}$\,.
$M_1$ is clearly a linear subspace and one can show that if $(f^{(k)})$
is a bounded sequence in $M_1$\,, converging pointwise to a
continuous function $f$\,, then
$\lVert\theta_f\rVert\le\limsup\lVert\theta_{f^{(k)}}\rVert$. In particular,
if  $(f^{(k)})$ is a Cauchy-sequence, then $f\in M_1$ and $f^{(k)}\to f$
in $M(A(G))$.

For $f=t_{jj'}^l$, with $l=-\frac12+i\lambda$ the evaluation of
$(f\odot\zeta\,,S)$ follows from \eqref{pi0} and the Corollary of
Lemma\;\ref{ebar}. Thus $f\in M_1$\,. This works in a similar way for the
coefficients of discrete series representations. 
Now observe that (using the formulas above) for fixed $m,n$ the function
$\lambda\mapsto t_{mn}^{-\frac12+i\lambda}\vert H$ \,($\R\to A(H)$) is
continuous. Then for $(m,n)$-radial $f\in A(G)$ one can approximate
the Fourier transform by finitely supported measures, giving an
approximation of $f$ (in $ M(A(G))$-norm) by linear combinations of
coefficients
$t_{mn}^l$\,. By the properties above this implies $f\in M_1$ and this
extends to arbitrary $f\in A(G)$ (and its norm closure in $M(A(G))$).
For general $f\in M(A(G))$ such that $f_{mn}\in A(H)$
for all $m,n\in2\Z$, one can use approximations (e.g. by Fejer sums) and
the properties of $M_1$ above to
see that $\Phi(f)$ belongs to the predual and
$\lVert\theta_f\rVert\le\lVert f\rVert_M$.

In the case of the $(n,m)$-radial functions $f'_{mn}=\chi_n*f*\chi_m$\,,
it follows easily (using $G=HK$\,,\ $f_{mn}\in B(H)$\,) that $f'_{mn}$ is
weakly almost periodic and
for $f$ $K$-finite, this implies that $f$ is weakly almost periodic.
By the results of \cite{Ve} it follows that $\lambda=\lim_{x\to\infty}f(x)$
exists and $f_0\in C_0(G)$\,.
For general $f\in M(A(G))$ this implies that $f_{mn}\in A(H)$ for
$(m,n)\ne(0,0)$ and there exists $\lambda\in\C$ such that
$(f-\lambda)_{00}=f_{00}-\lambda\in A(H)$\,. Then the
formulas involving $\Phi_1(f)$ and $\theta_{1f}$ follow first for
$K$-finite $f$\,, applying the previous results to $f_0$\,. Finally,
approximation gives the general case of the Corollary and the earlier
argument, based on Proposition 1 shows $f-\lambda\in C_0(G)$
\,(i.e., $\lambda=\lim_{x\to\infty}f(x)$\,)\vspace{1mm}.

As mentioned before we have restricted
to representations of $\PSL(2,\R)$ and this produces only $(m,n)$-radial
functions with $m,n$ even; the other representations of $\SL(2,\R)$ give
odd values for $m,n$ and this amounts to extend the definition of
$\overline{\Cal H}_p\,,\,\varphi,\dots$ to half-integer $j,s$\,.
\end{proof}
\newpage

\end{document}